\documentclass[12pt]{article}
\usepackage{amssymb}
\newcommand{\R}{\mathbb{R}}

\newcommand{\D}{\mathbb{D}}

\newcommand{\M}{\mathbb{M}}

\def\build#1_#2^#3{\mathrel{
\mathop{\kern 0pt#1}\limits_{#2}^{#3}}}

\def\cq{$\hfill \square$}

\def\t{{\cal T}}

\def\W{{\cal W}}
\def\w{{\rm w}}

\def\varep{\varepsilon}
\def\PP{\hbox{\bf P}}

\def\EE{\hbox{\bf E}}

\def\be{\begin{equation}}
\def\ee{\end{equation}}
\def\ba{\begin{eqnarray*}}
\def\ea{\end{eqnarray*}}
\def\ov{\overline}

\def\wt{\widetilde}

\def\proof{\vskip 3mm \noindent{\bf Proof:}\hskip10pt}

\newtheorem{theorem}{Theorem}[section]
\newtheorem{lemma}[theorem]{Lemma}
\newtheorem{proposition}[theorem]{Proposition}

\begin{document}

\title{ \bf On the re-rooting invariance property \\of L\'evy trees}
\author{Thomas {\sc Duquesne} and
Jean-Fran\c cois {\sc Le Gall}}

\vspace{2mm}
\date{\tiny\today} 

\maketitle

\begin{abstract}
We prove a strong form of the invariance under re-rooting of
the distribution of the continuous random trees called 
L\'evy trees. This extends previous results due to several authors.
\end{abstract}

\section{Introduction}

Continuous random trees have been studied extensively in the
last fifteen years and have found applications in various areas of
probability theory and combinatorics. The prototype of these random trees
is the CRT (Continuum Random Tree) which was introduced 
and studied by Aldous in a series of papers \cite{Al1,Al2,Al3}
in the early nineties. Amongst other things, Aldous proved that
the CRT is the scaling limit of a wide class of discrete random trees
(in particular it is the scaling limit of uniformly distributed planar trees with $n$ edges,
or of Cayley trees on $n$ vertices, when
$n$ tends to infinity) and that it can be coded by a normalized 
Brownian excursion. The latter fact makes it possible to derive a
number of explicit distributions for the CRT. 

The class of L\'evy trees, generalizing the CRT, was introduced and
discussed in the monograph \cite{DLG0}, and further studied in \cite{DLG}.
L\'evy trees are continuous analogues of Galton-Watson branching
trees and, in some sense, they are the only
continuous random trees that can arise as scaling limits of sequences of
Galton-Watson trees conditioned to be large: See Chapter 2
in \cite{DLG0} for a thorough discussion of these scaling limits.
Recently, Weill \cite{W} also proved that L\'evy trees are precisely those
continuous random trees that enjoy the regeneration property,
saying informally that the subtrees originating from a fixed level
in the tree are independent and distributed as the full tree.
Among L\'evy trees, the special class of stable trees,
which arise as scaling limits of sequences of conditioned Galton-Watson trees
with the same offspring distribution \cite{Duq}, is of
particular importance.

In the formalism of \cite{DLG}, or in \cite{EPW}, continuous random
trees are viewed as random variables taking values in the
space of all compact rooted real trees, which is itself a closed subset of the
space of all (isometry classes of) pointed compact metric spaces, which is equipped with the
Gromov-Hausdorff distance \cite{Gro}. Since the latter space is Polish,
this provides a convenient setting for the 
study of distributional properties of continuous random trees.

In a way analogous to the coding of discrete planar trees by
Dyck paths, one can code a compact real tree by a contour function
(see Theorem 2.1 in \cite{DLG}). Roughly speaking, with any
continuous function $g:[0,\sigma]\to\R_+$
such that $g(0)=g(\sigma)=0$ one associates a real tree $\t_g$
called the tree coded by $g$, which is formally obtained
as the quotient of the interval $[0,\sigma]$ for an
equivalence relation $\sim_g$ defined in terms of $g$ -- See Section 2 below. By definition
the root of $\t_g$ is the equivalence class of $0$, and the 
tree $\t_g$ is equipped with a ``uniform'' measure ${\bf m}$,
which is just the image of Lebesgue measure under the
canonical projection $[0,\sigma]\to [0,\sigma]\,/\!\sim_g=\t_g$. This construction 
applies to the CRT, which is the tree coded by a normalized 
Brownian excursion (thus $\sigma=1$ in that case), and more generally to L\'evy trees, for which
the coding function is the so-called height process, which is
a functional of a L\'evy process without negative jumps: See \cite{LGLJ1,DLG0,DLG}
and Section \ref{Not} below.

The main purpose of the present note is to discuss a remarkable
invariance property of L\'evy trees under re-rooting. 
In the case of the CRT, Aldous \cite{Al2} already observed that 
the law of the CRT is invariant under uniform re-rooting: If we pick 
a vertex of the CRT according to the uniform measure ${\bf m}$,
the CRT re-rooted at this vertex has the same distribution
as the CRT. As was noted by Aldous, this property corresponds
to the invariance of the law of the Brownian excursion
under a simple path transformation. An analogue of this property for L\'evy trees
was derived in Proposition 4.8 of \cite{DLG}. It turns out that
much more is true. Marckert and Mokkadem \cite{MaMo} (see also
Theorem 2.3 in \cite{LGW})
proved that,
for any fixed $s\in[0,1]$, the CRT re-rooted at the vertex $s$ (or rather the equivalence class
of $s$) has the same distribution as the CRT. The main result of the present work
(Theorem \ref{mainC})
shows that this strong form of the re-rooting invariance property
remains valid for L\'evy trees. In fact, we obtain the corresponding
result for the height process coding the L\'evy tree, which gives a more precise
statement than just the invariance of the distribution of the tree
(compare the first and the second assertion of Theorem \ref{mainC}).

In the case of the CRT, the invariance under re-rooting can be 
deduced from a similar property for approximating discrete random
trees. This approach was used both by  Aldous and by
Marckert and Mokkadem. In the case of L\'evy trees, this method seems
much harder to implement and we prefer to argue directly on the
coding function, using fine properties  of the underlying L\'evy
process. In a sense, our arguments are in the same spirit as the proof of
Theorem 2.3 in \cite{LGW}, but the latter paper used very 
specific properties of Brownian motion, which no longer hold
in our general setting.

Part of the motivation for the present work came from the recent paper
by Haas, Pitman and Winkel \cite{HPW}. This paper proves that,
among the continuous fragmentation trees, stable trees 
are the only ones whose distribution is invariant under uniform re-rooting.
Note that continuous random trees modelling the genealogy
of self-similar fragmentations were investigated by Haas and Miermont \cite{HaM}.
The present work gives more insight in the
probabilistic properties of the L\'evy trees, which are still the subject
of active research (see in particular the recent paper \cite{AD}).
The re-rooting invariance of the CRT turned out to play a
very important role in the study of certain conditionings considered
in \cite{LGW}, which are closely related to the continuous limit
of random planar maps \cite{MaMo,LGmap}. We hope that our
results about L\'evy trees will have similar applications 
to related models.

The paper is organized as follows. Section \ref{Not} recalls
the basic notation and the construction of the L\'evy tree
from the height process, and also shows how our main result
can be deduced from the technical Proposition \ref{mainT}.
At the end of Section 2, we briefly discuss some applications of
the invariance under re-rooting. Finally,
Section \ref{Proof} gives the proof of Proposition \ref{mainT}.

\noindent{\bf Acknowledgement.} We thank B\'en\'edicte Haas for
keeping us informed of her joint work with Jim Pitman
and Matthias Winkel.

\section{Notation and statement of the main results}
\label{Not}

We consider the general framework of 
\cite{DLG0,DLG}. Let $X=(X_t)_{t\geq 0}$ be a real L\'evy process
 without negative jumps, which starts from $0$ under the probability measure $\PP$. Without loss of generality,
 we may and will assume that $X$ is the canonical process on the Skorokhod
 space $\D(\R_+,\R)$. The 
Laplace exponent of $X$ is denoted by $\psi(\lambda)$:
$$\EE[\exp(-\lambda X_t)]=\exp(t\psi(\lambda))\;,\quad t,\lambda\geq 0.$$
We assume that $X$ has first moments and that $\EE[X_1]\leq 0$. Equivalently, this means
that we exclude the case where $X$ drifts to $+\infty$. Then $\psi$ can be written in
the form
$$ \psi (\lambda )= \alpha \lambda + \beta \lambda^2 + 
\int_{(0, \infty)} (e^{-\lambda r} -1+\lambda r)\,\pi (dr)  \; ,  $$
where $\alpha, \beta \geq 0$ and the L\'evy measure $\pi $ is a $\sigma $-finite measure on $(0, \infty)$ such that 
$\int_{(0, \infty)} (r\wedge r^2)\,\pi (dr) < \infty $. Let
$$I_t=\inf_{0\leq s\leq t}X_s\;,\quad t\geq 0\;,$$
denote the minimum process of $X$. Then the process
$X-I$ is a strong Markov process in $[0,\infty)$ and the point $0$ is regular and recurrent for this
Markov process (see Sections VI.1 and VII.1 in \cite{Be}).

We furthermore assume that
\begin{equation}
\label{extinct}
\int_1^\infty {du\over \psi(u)}<\infty.
\end{equation}
This implies that at least one of the following two conditions holds: 
\begin{equation}
\label{infvar}
\beta >0 \quad {\rm or} \quad \int_{(0, 1)} r\,\pi (dr)  = \infty .
\end{equation}
In particular (\cite{Be}, Corollary VII.5), the paths of $X$ are of infinite variation
almost surely, and the point $0$ is regular for $(0,\infty)$ with respect to $X$,
and thus also with respect to $X-I$. By a simple duality argument,
this implies that $\int_0^\infty {\bf 1}_{\{X_s=I_s\}} ds =0$.

By Theorem VII.1 in \cite{Be},  the continuous increasing process $(-I_t)_{t\geq 0}$
is a local time at $0$ for $X-I$. We denote by $N$ the associated It\^o excursion 
measure, so that $N$ is a $\sigma$-finite measure on the space 
$\D(\R_+,\R)$. We also denote by $\sigma:=\inf\{t>0:X_t=0\}$ the duration of
the excursion under $N$. 

Let us now introduce the height process. For every $0\leq s\leq t$, set
$$I^s_t=\inf_{s\leq r\leq t} X_r\;.$$
One can then prove (see Chapter 1 of \cite{DLG0}) that, for every $t\geq 0$, the limit
\begin{equation}
\label{Hlimit}
H_t:=\lim_{\varepsilon\to 0} \frac{1}{\varep}\int_0^t {\bf 1}_{\{I^s_t<X_s<I^s_t+\varep\}}\,ds
\end{equation}
exists in $\PP$-probability and in $N$-measure. The reason why it is useful to
consider the process $H$ also under the excursion measure $N$
is the fact that $\PP$ a.s. $H_t$ only depends on the excursion of $X-I$
that straddles time $t$: In the integral appearing in the
right-hand side of (\ref{Hlimit}), we can restrict our attention to values
of $s$ belonging to the excursion interval of $X-I$ away from $0$
that straddles $t$. See the discussion in Section 3.2 of \cite{DLG}.

Thanks to
our assumption (\ref{extinct}), the process $(H_t)_{t\geq 0}$ has a continuous
modification under $\PP$ and under $N$: This means that we can redefine
the process $(H_t)_{t\geq 0}$ on the Skorokhod space in such a way that
all its sample paths are continuous and the limit (\ref{Hlimit}) still holds for
every $t\geq 0$. Notice that $H_0=0$, $\PP$ a.s. and $N$ a.e., and that
$H_t=0$ for every $t\geq \sigma$, $N$ a.e. The time-reversal invariance property
of $H$ states that the two processes $(H_t)_{t\geq 0}$ and 
$(H_{(\sigma-t)_+})_{t\geq 0}$ have the same distribution under $N$
(see Corollary 3.1.6 in \cite{DLG0}).

The $\psi$-{\it L\'evy tree} is by definition the rooted real tree $\t_H$ coded by the 
continuous function $(H_t)_{t\geq 0}$ under the measure $N$, in the sense
of Theorem 2.1 in \cite{DLG}. This means that $\t_H=[0,\sigma]\,/\!\sim_H$, where the
random equivalence relation $\sim_H$ is defined on $[0,\sigma]$ by
$$s\sim_H t\hbox{\quad iff\quad} H_s=H_t=\min_{s\wedge t\leq r\leq s\vee t} H_r\;.$$
The canonical projection from $[0,\sigma]$ onto $\t_H$ is
denoted by $p_H$, and the distance on $\t_H$ is given by
$$d_H(p_H(s),p_H(t))=d_H(s,t)= H_s+H_t -2 \min_{s\wedge t\leq r\leq s\vee t} H_r\;$$
(notice that this only depends on $p_H(s)$ and $p_H(t)$, and not on the particular choice
of the representatives $s$ and $t$). By definition the tree $\t_H$
is rooted at $\rho=p_H(0)=p_H(\sigma)$.

In order to discuss re-rooting, let us introduce the following notation. Fix $s\in[0,\sigma]$
and set
$$H^{[s]}_t=\left\{
\begin{array}{ll}
d_H(s,s+t)\quad&\hbox{if }0\leq t<\sigma -s\\
d_H(s,s+t-\sigma)\quad&\hbox{if }\sigma-s\leq t\leq \sigma
\end{array}
\right.
$$
and $H^{[s]}_t=0$ if $t>\sigma$. By Lemma 2.2 in \cite{DLG}, the tree 
$\t_{H^{[s]}}$ is then canonically identified with the tree $\t_H$ re-rooted
at the vertex $p_H(s)$. 

We can now state the key technical
proposition that will lead to our main theorem. We denote by $C_0(\R_+,\R_+)$ the space of all
continuous functions with compact support from $\R_+$ into $\R_+$, so that
$(H_t)_{t\geq 0}$, or $(H^{[s]}_t)_{t\geq 0}$, can be viewed under $N$ as a random
element of $C_0(\R_+,\R_+)$.

\begin{proposition}
\label{mainT}
For every nonnegative measurable function $F$
on $\R_+\times C_0(\R_+,\R_+)$, and every nonnegative measurable function
$g$ from $\R_+$ into $\R_+$,
$$N\Big(g(\sigma)\int_0^\sigma ds\,F(s,H^{[s]})\Big)=N\Big(g(\sigma)\int_0^\sigma ds\,F(s,H)\Big).$$
\end{proposition}

The proof of Proposition
\ref{mainT} is given in the next section, but we will immediately
deduce our main theorem from this proposition. Denote by $\kappa(da)$ the ``law'' of $\sigma$ under $N$. By a standard 
disintegration theorem, we can find a measurable collection $(N^{(a)})_{a\in(0,\infty)}$
of probability measures on $\D(\R_+,\R)$ such that:
\begin{itemize}
\item[(a)] $N=\int_0^\infty \kappa(da)\,N^{(a)}\;$;
\item[(b)] for every $a\in(0,\infty)$, $N^{(a)}$ is supported on $\{\sigma = a\}$.
\end{itemize}
In the stable case where $\psi(u) = c \,u^\gamma$ for 
some $c>0$ and $\gamma\in(1,2]$, the measures $N^{(a)}$ can be chosen in such a way that
$N^{(a)}$ is the law of $(a^{1/\gamma}X_{t/a})_{t\geq 0}$ under $N^{(1)}$. 
A scaling argument then shows that the continuous modification of $H$ can be chosen so that
(\ref{Hlimit}) holds in $N^{(a)}$-probability for every $a>0$. So we may, and will, assume
that the latter properties hold in the stable case. 

By definition, the {\it stable tree}
is the tree $\t_H$ under the probability measure $N^{(1)}$, when
$\psi(u)=u^\gamma$ for some
$\gamma\in (1,2]$. The CRT is the case  $\gamma=2$. 

\begin{theorem}
\label{mainC}
The following properties hold for $\kappa$-almost every $a>0$. For every 
$s_0\in[0,a):$
\begin{description}
\item{\rm(i)}
The processes $(H^{[s_0]}_t)_{t\geq 0}$ and $(H_t)_{t\geq 0}$ 
have the same distribution under $N^{(a)}$. 
\item{\rm(ii)} The law under
$N^{(a)}$ of the tree $\t_H$ re-rooted at $p_H(s_0)$ coincides with the law
of $\t_H$.
\end{description}
In the stable case, the preceding properties hold for every $a>0$. 
\end{theorem}

\noindent{\bf Remark.} Theorem \ref{mainC} generalizes property (20) in \cite{Al2} and
Proposition 4.9 in \cite{MaMo}, which are both concerned with Aldous' CRT, as well
as Theorem 11(i) in \cite{HPW}, which deals with uniform re-rooting of the stable tree. 
Theorem \ref{mainC}
also strengthens Proposition 4.8 of \cite{DLG}, which considers only uniform
re-rooting and property (ii).

\medskip
\proof Using properties (a) and (b), we immediately get that, for every fixed function $F$
satisfying the properties in Proposition \ref{mainT}, the equality
$$N^{(a)}\Big(\int_0^a ds\,F(s,H^{[s]})\Big)=N^{(a)}\Big(\int_0^a ds\,F(s,H)\Big)$$
holds $\kappa(da)$ a.e. A simple separability argument implies that the preceding
identity holds simultaneously for all choices of $F$, except for values of $a$
belonging to a set
of zero $\kappa$-measure. We apply this with $F(s,H)=\varep^{-1}{\bf 1}_{[s_0,s_0+\varep]}(s)G(H)$
where $G$ is continuous and bounded on $C_0(\R_+,\R_+)$. Letting
$\varep$ go to $0$, and noting that $H^{[s]}$ depends continuously on $s$, it follows that
for $\kappa$-almost every $a>0$, we have for every such function $G$,
for every $s_0\in[0,a)$,
$$N^{(a)}(G(H^{[s_0]}))= N^{(a)}(G(H)).$$
Property (i) readily follows. Property (ii) is then a consequence of the
fact that $\t_{H^{[s_0]}}$ is isometric to the tree $\t_H$ re-rooted at $p_H(s_0)$
(\cite{DLG}, Lemma 2.2). 

In the stable case, a scaling argument shows that the properties of the theorem hold for
every $a>0$ as soon as they hold for one value of $a>0$. \cq

\medskip
Let us briefly comment on applications of the invariance under re-rooting. Aldous \cite{Al2,Al3}
observed that a convenient way to describe the distribution of a continuous random tree
is via its finite-dimensional marginals. Let $\t$ be a rooted continuous random tree (for instance
the stable tree) and suppose that we are given a probability distribution ${\bf m}$ on $\t$,
which is called the mass measure (in the case of the stable tree, ${\bf m}$ is the
image of Lebesgue measure over $[0,1]$ under the projection $p_H$). Let $p\geq 1$ be an integer 
and suppose that $V_1,\ldots,V_p$ are $p$ vertices chosen independently at random on $\t$
according to the mass measure ${\bf m}$. The subtree $\t(V_1,\ldots,V_p)$
 spanned by $V_1,\ldots,V_p$ is the union of the line 
segments beween the root and the vertices $V_1,\ldots,V_p$. It is a finite real tree, 
meaning that it consists of the union of a finite number of segments, with 
$p+1$ labeled vertices: the root is labeled $0$ and each vertex $V_k$
is labeled $k$. The invariance under uniform re-rooting implies that the 
distribution of this labeled tree is invariant under every permutation of
the labels $0,1,\ldots,p$. For instance, in the case when $\t$ is the stable tree
and $p=2$, we get that if $U$ and $U'$ are two independent random variables
uniformly distributed over $[0,1]$,
and if $m_H(U,U')=\min_{U\wedge U'\leq s\leq U\vee U'} H_s$, the law under $N^{(1)}$ of the triplet
$$\Big(H_{U}-m_H(U,U')\;,\,
H_{U'}-m_H(U,U')\;,\,m_H(U,U')\Big)$$
is exchangeable. Even in the case of the CRT, where $H$ is the normalized Brownian
excursion, it is not so easy to give a direct derivation of this property.

Other applications of the invariance under re-rooting are concerned with tree-indexed
processes (see in particular \cite{LGW}). In order to provide a simple example,
consider again the stable tree $\t_H$,
under the probability measure $N^{(1)}$ as above, and its mass measure ${\bf m}$. Let $(Z_a,a\in\t_H)$ be Brownian motion indexed by $\t_H$, which may
defined as follows. Conditionally given $\t_H$, this is the centered Gaussian process
such that $Z_\rho=0$ and $E[(Z_a-Z_b)^2]=d_H(a,b)$ for every $a,b\in\t_H$. The 
occupation measure ${\cal I}$ of $Z$ is the random measure on
$\R$ defined by:
$$\langle {\cal I},g\rangle = \int {\bf m}(da) \,g(Z_a).$$
In the particular case of the CRT, the measure $\cal I$ is known as
(one-dimensional) ISE -- See Aldous \cite{Al4}. As a consequence of 
the invariance under uniform re-rooting, one easily checks that the
quantity ${\cal I}(]0,\infty[)$, which represents the total mass to the right
of the origin, is uniformly distributed over $[0,1]$. This had already been
observed in \cite{Al4} in the case of the CRT.

\section{Proof of the main proposition}
\label{Proof}

In this section we prove Proposition \ref{mainT}. We start by recalling a lemma
from \cite{DLG} that plays a key role in our proof. We first need to introduce some
notation. It will be convenient to use the notion of a finite path. A (one-dimensional)
finite path is just a continuous mapping $\w:[0,\zeta]\longrightarrow \R$, where the
number $\zeta=\zeta_{(\w)}\geq 0$ is called the lifetime of $\w$. The space $\W$
of all finite paths is equipped with the distance $d$ defined by
$d(\w,\w')=\|\w-\w'\| + |\zeta_{(\w)}-\zeta_{(\w')}|$, where $\|\w-\w'\|
=\sup_{t\geq 0}|\w(t\wedge\zeta_{(\w)})-\w'(t\wedge\zeta_{(\w')})|$.

Let $M_f$ denote the space of 	all finite measures on $\R_+$. For every $\mu\in M_f$, let
${\rm supp}(\mu)$ denote the topological support of $\mu$ and set $S(\mu)=\sup({\rm supp}(\mu))
\in[0,+\infty]$. We let $M_f^*$ be the subset of $M_f$ consisting of
all measures $\mu$ such that $S(\mu)<\infty$ and ${\rm supp}(\mu)=[0,S(\mu)]$. By convention
the measure $\mu=0$ belongs to $M_f^*$ and $S(0)=0$. 

Let $\psi^*(u)=\psi(u)-\alpha\,u$, and let 
$(U^1,U^2)$ be a two-dimensional subordinator with Laplace functional
$$E[\exp(-\lambda U^1_t-\lambda'U^2_t)]=\exp\Big(-{\psi^*(\lambda)-\psi^*(\lambda')\over \lambda-\lambda'}
\Big).$$
Note that $U^1$ and $U^2$ have the same distribution and are indeed
subordinators with drift $\beta$ and L\'evy measure $\pi([x,\infty))dx$. 
For every $a> 0$, we let $\M_a$ be the
probability measure on
$(M^*_f)^2$ which is the distribution of $({\bf 1}_{[0,a]}(t)\,dU^1_t,{\bf 1}_{[0,a]}(t)\,dU^2_t)$.
The fact that this defines a distribution on $(M^*_f)^2$ follows from (\ref{infvar}).
Moreover $\M_a(d\mu\,d\nu)$ a.s., we have $S(\mu)=S(\nu)=a$. 

Let $\mu\in M^*_f$ and denote by $|\mu|=\mu([0,\infty))$ the total mass of $\mu$. For every
$r\in[0,|\mu|]$, we denote by $k_r\mu$ the unique element
of $M_f^*$ such that
$$k_r\mu([0,x])=\mu([0,x])\wedge (|\mu|-r)$$
for every $x\geq 0$. We then set, for every $0\leq t\leq T_{|\mu|}:=\inf\{s\geq 0:X_s=-|\mu|\}$,
$$H^\mu_t=S(k_{-I_t}\mu)+H_t.$$
We view $(H^\mu_t)_{0\leq t\leq T_{|\mu|}}$ as a random element of $\W$,
with lifetime $\zeta_{(H^\mu)}=T_{|\mu|}$, and we let $Q_\mu(d\w)$
be the probability measure on $\W$ which is the distribution of $H^\mu$
under $\PP$. 
Notice that $\w(0)=S(\mu)$ and $\w(\zeta_{(\w)})=0$, $Q_\mu(d\w)$ a.s.

Finally, under the excursion measure $N$, we set for every 
$s\in[0,\sigma]$,
$$\begin{array}{ll}
H^{+,s}_t=H_{s+t}\;,\quad&0\leq t\leq \sigma-s\;,\\
H^{-,s}_t=H_{s-t}\;,\quad&0\leq t\leq s\;,
\end{array}
$$
and we view $H^{+,s}$ and $H^{-,s}$ as
random elements of $\W$ with respective lifetimes
$\sigma - s$ and $s$.

The following result is Lemma 3.4 in \cite{DLG}.

\begin{lemma}
\label{keytech}
For any nonnegative measurable function $\Phi$ on $\W^2$,
$$N\Big(\int_0^\zeta ds\,\Phi(H^{+,s},H^{-,s})\Big)
=\int_0^\infty da\,e^{-\alpha a}\int \M_a(d\mu d\nu)\int \!\!\int Q_\mu(d\w)Q_\nu(d\w')\Phi(\w,\w').
$$
\end{lemma}

Before stating the next lemma we introduce two simple transformations
of finite paths. Let $\w\in\W$. We define two other finite paths
$\ov \w$ and $\wt \w$ both having the same
lifetime as $\w$, by setting
$$\begin{array}{ll}
\ov\w(t)=\w(\zeta_{(\w)}-t)\;,\quad&0\leq t\leq \zeta_{(\w)}\;,\\
\noalign{\smallskip}
\wt\w(t)=\w(0)+\w(t)-2
\displaystyle{\min_{0\leq s\leq t}\w(s)}\;,\quad&0\leq t\leq \zeta_{(\w)}\;.
\end{array}
$$
If $\mu\in M_f^*$, we also denote by $\ov \mu$ the ``time-reversed'' measure
defined as the image of $\mu$ under the mapping $t\to S(\mu)-t$.

\begin{lemma}
\label{key2}
Let $\mu\in M_f^*$. The law of $\wt\w$ under $Q_\mu(d\w)$ coincides with the
law of $\ov \w$ under $Q_{\ov\mu}(d\w)$. 
\end{lemma}

\proof To simplify notation, we write $\tau=T_{|\mu|}$ in this proof.
We argue under the probability measure $\PP$. Since the
function $t\to S(k_{-I_t}\mu)$ is nonincreasing and can only
decrease when $X_t=I_t$, which forces $H_t=0$, it is easy to verify that,
a.s. for every $t\in[0,\tau]$,
$$\min_{0\leq r\leq t} H^\mu_r=S(k_{-I_t}\mu)$$
(consider the last time $r$ before $t$ when $X_r=I_r$). It follows that
\begin{equation}
\label{key2-1}
\wt H^\mu_t=H^\mu_0 + H^\mu_t -2 \min_{0\leq r\leq t} H^\mu_r = S(\mu) +H_t-S(k_{-I_t}\mu),
\end{equation}
for every $t\in[0,\tau]$, a.s. 

Now observe that the random function  $(H_t, t\in[0,\tau])$ is obtained by
concatenating the values of $H$ over all excursion intervals of $X-I$
away from $0$ between times $0$ and $\tau$. Since 
$-I$ is a local time at $0$
for the process $X-I$, and $\tau=\inf\{t\geq 0: -I_t=|\mu|\}$, the excursions
of $X-I$ between times $0$ and $\tau$ form a Poisson
point process with intensity $|\mu|\,N$. Furthermore the values of 
$H$ over each excursion interval only depend on the corresponding excursion of $X-I$, and
we know that the law of $H$ under $N$ is invariant under time-reversal.
By putting the previous facts together, and using standard 
arguments of excursion theory, we obtain that the two pairs
of processes
$$(H_{\tau-t},|\mu|+I_{\tau-t})_{0\leq t\leq \tau}\quad
\hbox{and}\quad (H_t,-I_t)_{0\leq t\leq \tau}$$
have the same distribution under $\PP$. 

From this identity in distribution, we deduce that, again under the probability measure $\PP$,
\begin{equation}
\label{key2-2}
(S(\mu)+H_{t}-S(k_{-I_{t}}\mu))_{0\leq t\leq \tau}
\build{=}_{}^{\rm(d)}
(S(\mu)+H_{\tau -t}-S(k_{|\mu|+I_{\tau-t}}\mu))_{0\leq t\leq \tau}.
\end{equation}
However, the elementary fact $S(\mu)-S(k_{|\mu|-r}\mu)=S(k_r\ov\mu)$, which is valid
for every $r\in[0,|\mu|]$, implies that, for every $t\in[0,\tau]$,
\begin{equation}
\label{key2-3}
S(\mu)+H_t-S(k_{|\mu|+I_t}\mu)=H_t+S(k_{-I_t}\ov\mu)=H^{\ov\mu}_t.
\end{equation}
The lemma now follows from (\ref{key2-1}), (\ref{key2-2}), (\ref{key2-3}) and the
definition of the measures $Q_{\mu}$. \cq

We now turn to the proof of Proposition \ref{mainT}. We first observe that, 
with the notation introduced before Lemma \ref{key2},
we have $N$ a.e. for every $s\in(0,\sigma)$,
$$
\begin{array}{ll}
\wt H^{+,s}_t = H^{[s]}_t\;,\quad &0\leq t\leq \sigma - s\;,\\
\wt H^{-,s}_t=H^{[s]}_{\sigma-t}\;,\quad &0\leq t\leq s\;.
\end{array}
$$
Let $\Psi$ be a nonnegative measurable function on $\W^2$. Using 
Lemma \ref{keytech} in the first equality and Lemma \ref{key2}
in the second one, we have
\newpage
\begin{eqnarray*}
&&N\Big(\int_0^\sigma ds\,\Psi\Big((H^{[s]}_t)_{0\leq t\leq \sigma -s},
(H^{[s]}_{\sigma-t})_{0\leq t\leq s}\Big)\Big)\\
&&\qquad=\int_0^\infty da\,e^{-\alpha a}\int \M_a(d\mu d\nu)
\int\!\!\int Q_\mu(d\w)Q_\nu(d\w')\,\Psi(\wt\w,\wt\w')\\
&&\qquad= \int_0^\infty da\,e^{-\alpha a}\int \M_a(d\mu d\nu)
\int\!\!\int Q_{\ov\mu}(d\w)Q_{\ov\nu}(d\w')\,\Psi(\ov\w,\ov\w')\\
&&\qquad= \int_0^\infty da\,e^{-\alpha a}\int \M_a(d\mu d\nu)
\int\!\!\int Q_\mu(d\w)Q_\nu(d\w')\,\Psi(\ov\w',\ov\w)\\
&&\qquad=N\Big(\int_0^\sigma ds\,\Psi(\ov H^{-,s},\ov H^{+,s})\Big)\\
&&\qquad=N\Big(\int_0^\sigma ds\,\Psi
\Big((H_t)_{0\leq t\leq s},(H_{\sigma-t})_{0\leq t\leq \sigma-s}\Big)\Big)
\end{eqnarray*}
In the third equality, we use the fact that the probability measure 
$\M_a(d\mu d\nu)$ is symmetric and invariant under the mapping
$(\mu,\nu)\to (\ov\mu,\ov \nu)$. The fourth equality is 
Lemma \ref{keytech} again, and the last one
is an immediate consequence of the definitions.

Now note that the triplet $(\sigma,\sigma-s,H)$ 
can be written as a measurable functional $\Gamma$ of the pair
$((H_t)_{0\leq t\leq s},(H_{\sigma-t})_{0\leq t\leq \sigma-s})$, and that
the triplet $(\sigma,s,H^{[s]})$ is then the same measurable function
of the pair $((H^{[s]}_t)_{0\leq t\leq \sigma -s},
(H^{[s]}_{\sigma-t})_{0\leq t\leq s}))$. Therefore, the preceding
calculation implies that, for $F$ and $g$ as in the
statement of the proposition,
$$N\Big(g(\sigma)\int_0^\sigma ds\,F(s,H^{[s]})\Big)=
N\Big(g(\sigma)\int_0^\sigma ds\,F(\sigma-s,H)\Big)=
N\Big(g(\sigma)\int_0^\sigma ds\,F(s,H)\Big).$$
This completes the proof. \cq

\medskip

{\small 
\begin{tabular}{l}
Thomas DUQUESNE \\
Laboratoire de Probabilit\'es et Mod\`eles Al\'eatoires\\
Universit\'e Paris 6, 16 rue Clisson, 75013 PARIS, FRANCE\\
e-mail: thomas.duquesne@upmc.fr
\end{tabular}

\begin{tabular}{l}
Jean-Fran\c cois LE GALL\\ 
Math\'ematiques, Universit\'e Paris-Sud\\
Centre d'Orsay, 91405 ORSAY Cedex, FRANCE \\
e-mail: jean-francois.legall@math.u-psud.fr
\end{tabular}
}

\end{document}